\documentclass[10pt]{amsart}
\usepackage{amsmath,amsfonts,amssymb,latexsym,amsthm}
\usepackage{epsfig}
\usepackage{hyperref} 
\hypersetup{colorlinks} 
\usepackage[dvipsnames]{xcolor}
\usepackage{changepage}

\usepackage{dsfont}
\usepackage{mathtools}
\usepackage{fancyhdr}
\usepackage{enumitem}
\usepackage{pgfplots}
\usepackage{url}
\usepackage{tcolorbox}
\usetikzlibrary{graphs}
\usetikzlibrary{graphs.standard}
\usetikzlibrary{arrows}
\usepackage{xcolor}
\usetikzlibrary{arrows.meta}
\usetikzlibrary{graphs}
\usetikzlibrary{graphs.standard}
\usetikzlibrary{arrows}
\usetikzlibrary{cd}

\usepackage{enumitem}

\numberwithin{equation}{section}

\newtheorem{thm}{Theorem}[section]

\newtheorem{conv}[thm]{Convention}

\newtheorem{theorem}[thm]{Theorem}

\makeatletter
\newcommand{\eqnum}{\leavevmode\hfill\refstepcounter{equation}\textup{\tagform@{\theequation}}}\makeatother

\title[Complexity of quiver mutation equivalence]{Complexity of quiver mutation equivalence}

\author{David Soukup}
\thanks{\thinspace ${\hspace{-.45ex}}^\star$Department of Mathematics,
UCLA, Los Angeles, CA~90095.
\hskip.06cm
Email:
\hskip.06cm
\texttt{soukup@math.ucla.edu}}

\newcommand{\Tow}{\mathrm{Tow}}

\def\<{\langle}
\def\>{\rangle}

\def\0{{\mathbf 0}}

\def\.{\hskip.06cm}

\usepackage[colorinlistoftodos]{todonotes}

\RequirePackage{cleveref}
\usepackage{hypcap}
\hypersetup{colorlinks=true, citecolor=darkblue, linkcolor=darkblue}
\definecolor{darkblue}{rgb}{0.0,0,0.7}
\newcommand{\darkblue}{\color{darkblue}}

\definecolor{darkred}{rgb}{0.68,0,0}

\definecolor{darkgreen}{rgb}{0,.38,0}

\newcommand{\defn}[1]{\emph{\darkblue #1}}

\begin{document}

\maketitle

\begin{abstract}
We prove \textsf{NP-hardness} results for determining whether quivers are mutation equivalent to quivers with given properties. Specifically, determining whether a quiver is mutation-equivalent to a quiver with exactly $k$ arrows between any two of its vertices is \textsf{NP-hard}. Also, determining whether a quiver is mutation equivalent to a quiver with no edges between frozen vertices is  \textsf{strongly NP-hard}. Finally, we present a characterization of mutation classes of quivers with two mutable vertices.
\end{abstract}

\section{Introduction}
Quivers and their mutations (defined in Section \ref{ss:definitions}) were introduced by Fomin and Zelevinsky in \cite{FZ1} and \cite{FZ2} in the context of cluster algebras. They are widely used in algebraic combinatorics (see \cite{Kel} for a survey or \cite{FWZ} for an introductory book.). However, many combinatorial questions about these objects remain unresolved.

The question of whether a given quiver $Q$ is equivalent to only finitely many other quivers was addressed in \cite{FST}, where a list of such quivers is given. Also, Fomin and Neville show in \cite{FN23+} that there are long cycles in the graph of quivers. Recently, Fomin asked in \cite{Fom22} for algorithmic solutions to the following questions:

\begin{itemize}
\item[1.] Given quivers $Q_1$ and $Q_2$, determine whether  $Q_1$ and $Q_2$ are mutation equivalent.
\item[2.] Given a quiver $Q$ and a nonnegative integer $k$, determine whether there exists a quiver $Q' \in [Q]$ such that $Q'$ has two vertices with exactly $k$ arrows between them.
\end{itemize}

However, Fomin also proposed that these problems may be computationally difficult or even undecidable:

\smallskip
\begin{adjustwidth}{30pt}{30pt}
\begin{center}

\textit{``We don't have any algorithm that would detect if two quivers are mutation equivalent or not ... of course it would be absurd if this were algorithmically undecidable - there must be an algorithm - well, who knows? Maybe not." \qquad \qquad }

\vspace{1mm}

\qquad \qquad --Sergey Fomin, \cite[approx. 19:00]{Fom22}, May 16, 2022

\end{center}
\end{adjustwidth}

\smallskip
\smallskip

Formally, Fomin's problem asks whether these questions are decidable for general quivers (cf. Problem 2.6.14 and Remark 4.1.13 in \cite{FWZ}). We approach the problem from both ends.
First, we present a couple of \textsf{NP-hardness} results about the second of Fomin's questions. Past results such as \cite{BFG} have shown that certain determinants are preserved by quiver mutation. Since determinants can be computed in polynomial time, however, these results show that it is unlikely that a determinantal formula can capture everything that is going on in a quiver. Next, we will show that quivers with only two mutable vertices can only have a very limited set of equivalent quivers, and derive asymptotics of  the quivers in such mutation classes.

Note that we use a slight generalization of quivers in which we do not ignore edges between frozen vertices. Such quivers have been studied in the literature (in e.g. \cite{Pre20}).

\subsection{Hardness results}

We begin by stating our main results. Both concern complexity of questions related to quiver mutation equivalence, specifically Fomin's second question. We note that these results do assume that frozen vertices and arrows between them are permitted in quivers.

\smallskip

\begin{theorem}[\textsf{NP-hardness}] \label{t:nphard}
Let $Q$ be a quiver, and let $k > 1$ be an integer. The following problem is \textup{\textsf{NP-hard}}: Determine whether there exists a quiver $Q'$ which is mutation equivalent to $Q$ such that $Q'$ contains two vertices with exactly $k$ arrows between them.
\end{theorem}

\smallskip

In the context of quivers, it is natural to be interested in  \textsf{strong NP-hardness}. In ordinary \textsf{NP-hardness}, the inputs to the problem are assumed to be in binary. Specifically, when there are $k$ arrows between two vertices in a quiver, this is assumed to take  $\log_2 k$ bits of input. However, the arrows in a quiver may each carry algebraic information and thus have independent meaning. When inputs to a decision problem are given in unary instead of binary, then the corresponding notion is \textsf{strong NP-hardness}. Problems such as \textsc{Knapsack} or \textsc{Subset Sum} do not meet this stronger criterion. See \cite{GJ1} for background on this topic.

Let an arrow in a quiver be \defn{icebound} if it goes between two frozen vertices.

\smallskip

\begin{thm}[\textsf{Strong NP-hardness}]\label{t:strongnphard}
Let $Q$ be a quiver. The following problem is \textup{\textsf{strongly NP-hard}}: Determine whether there exists a sequence of mutations which takes $Q$ to a quiver with no icebound arrows.
\end{thm}
  
  See Section \ref{ss:pnp} for implications of these results.

\subsection{Asymptotic result}

Let us limit the number of vertices at which we are allowed to mutate the quiver, then the set of mutation-equivalent quivers becomes quite limited. If there is only one mutable vertex, then, since mutation is an involution, there can only be two quivers in a mutation class.

Our theorem describes the mutation classes of quivers with exactly two mutable vertices. Again, since mutations are involutions, the only way to get new quivers is to alternate mutating at the two vertices. 

\begin{thm}  \label{t:alternate}
Let $Q$ be a quiver with exactly two mutable vertices called $C$ and $D$. Define $\alpha$ to be the number of arrows between $C$ and $D$. Then: 

\begin{itemize}
\item[]If $\alpha = 0$, we have $\big| [Q] \big| \leq 4$.
\item[] If $\alpha = 1$, we have $\big| [Q] \big| \leq 10$.
\item[] If $\alpha = 2$, then in any nontrivial case the number of arrows in $(\mu_D\mu_C)^n(Q)$ grows linearly.
\item[] If $\alpha \geq 3$, then in any nontrivial case the number of arrows in $(\mu_D\mu_C)^n(Q)$ grows exponentially.
\end{itemize}

Furthermore, if $\alpha \geq 2$, let $\delta_{I, J}(n)$ be the number of arrows between vertices $I$ and $J$ in $(\mu_D\mu_C)^n(Q)$. For any vertex $A \neq C, D$ we have:
\[
\lim_{n \to \infty }~\frac{\delta_{A, C}(n)}{\delta_{A, D}(n)} = \frac{1}{2}\left(\alpha + \sqrt{\alpha^2 - 4} \right)
\]

\end{thm}

See \ref{ss:gadgets} for possible extensions of this result.

\subsection{Structure of the paper.} We will proceed as follows. In Section \ref{ss:definitions} we begin with notation, definitions, and examples. Next we prove our theorem in Section \ref{s:proofs}. We conclude with final remarks in Section \ref{s:finalremarks}.

\section{Notation, Definitions, and Examples} \label{ss:definitions}

\subsection{Basic definitions}
For positive integers $n$, define $[n]$ to be the set $\{1, 2, \dots, n\}$. Also, let $\mathbb{N}$ be the set $\{0, 1, 2, \cdots\}$.

\subsection{Quivers}
A \defn{quiver} is a directed multigraph with no loops or 2-cycles, the edges of which are called \defn{arrows}. We will indicate multiple arrows between vertices by labeling edges with numbers. For example, the following graph is a quiver on five vertices with eight arrows:

\smallskip

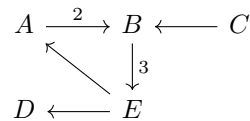
\begin{center}

\begin{tikzcd}
A \ar[r, "2"]& B\ar[d, "3"]  & C \ar[l]\\
D & E \ar[ul] \ar[l]
\end{tikzcd}

\end{center}

\smallskip

\subsection{Quiver mutation} In a quiver, we assign a subset of the vertices to be \defn{mutable}; the remaining vertices are \defn{frozen}. While Fomin and Zelevinsky's original definition ignored any arrows between frozen vertices, we will follow \cite{Pre20} and allow them. To each mutable vertex in the quiver we associate an operation called \defn{mutation}. For a vertex $X$, mutation at $X$, denoted by $\mu_X$, proceeds in the following three steps:

\smallskip

\begin{itemize}
\item[1.] for every two step path $Y \to X \to Z$, add an arrow from $Y$ to $Z$,
\item[2.] reverse the direction of every arrow incident to $X$,
\item[3.] remove 2-cycles one by one.
\end{itemize}

\smallskip

For example, applying the mutation $\mu_B$ will turn the quiver on the left into the quiver on the right and vice versa in the picture below:

\smallskip
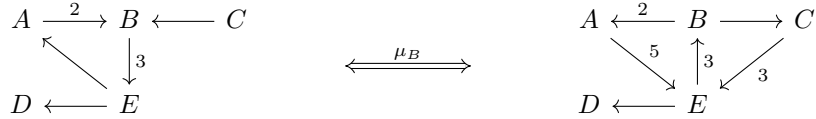
\begin{center}

\begin{tikzcd}
A \ar[r, "2"]& B\ar[d, "3"]  & C \ar[l]		
&\ar[rr, Leftrightarrow, shift right=7, "\mu_B"]&&\ & A \ar[dr, "5"]& B \ar[l, "2"']  \ar[r]  & C \ar[dl, "3"]\\
D & E \ar[ul] \ar[l]  &	
&&&& D & E  \ar[u, "3"'] \ar[l]   &
\end{tikzcd}

\end{center}

\smallskip

It is easily seen that every mutation is an involution. That is, $\mu_X(\mu_X(Q)) = Q$ for every quiver $Q$ with vertex $X$. It is also easy to see that mutations at nonadjacent vertices commute. Two quivers are said to be \defn{mutation equivalent} if one can be obtained from the other by a finite sequence of mutations.  Mutation equivalence is an equivalence relation, so we can define the \defn{mutation class} of a quiver $Q$, denoted $[Q]$, to be the equivalence class of $Q$ under this relation.

\smallskip

\section{Proofs} \label{s:proofs}

\subsection*{Proof of Theorem \ref{t:nphard}}

We reduce the problem to \textsc{Subset Sum}, which is defined as follows:

\smallskip

\begin{tabular}{l l}
\multicolumn{2}{l}{\textsc{\large{Subset Sum } } } \\
\textbf{Input:}  & $X \subset \mathbb{N}$ a finite set, and $k \in \mathbb{N}$.							\\ 
\textbf{Decide:} \quad & $\exists A \subseteq X$ ~ such that ~ $\sum_{a \in A}\  a = k$?

\end{tabular}

\smallskip

This problem is is known to be \textsf{NP-hard} (see e.g. \cite[\S A3.2]{GJ}). Let $X = \{x_1, \dots, x_n\}$ be a set of positive integers, and let $k > 1$ be   another integer. Let $Q$ be the following quiver:

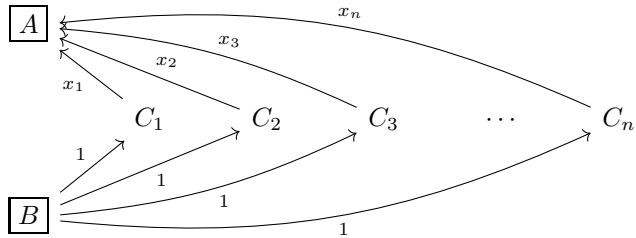
\vspace{5mm}

\begin{center}
\begin{tikzcd}
\boxed{A} \\
& C_1 \ar[ul, "x_1"] & C_2 \ar[ull, "x_2"'] & C_3 \ar[ulll, "x_3"', bend right = 10]& \cdots &  C_n \ar[ulllll,  "x_n"', bend right = 15]\\
\boxed{B} \ar[ur, "1"]  \ar[urr, "1"'] \ar[urrr, "1"', bend right = 10] \ar[urrrrr,  "1"', bend right = 15]
\end{tikzcd}
\end{center}

\vspace{5mm}

For each $i \in [n]$, let $\mu_i$ be $\mu_{C_i}$. Suppose we apply the sequence of mutations $\mu = \mu_{{i_1}} \cdots \mu_{{i_k}}$. Define $Y \subseteq [n]$ by
\[
Y = \{ j \in [n]\ : \ \mu_{C_j} \text{ is used an odd number of times} \}
\]
Then, for each $j \in [n]$ let 
\[
\epsilon_j = \begin{cases} 1 &\text{if } j \notin Y \\
-1 & \text{if } j \in Y
\end{cases}
\]
An easy induction shows that $\mu(Q)$ is given by

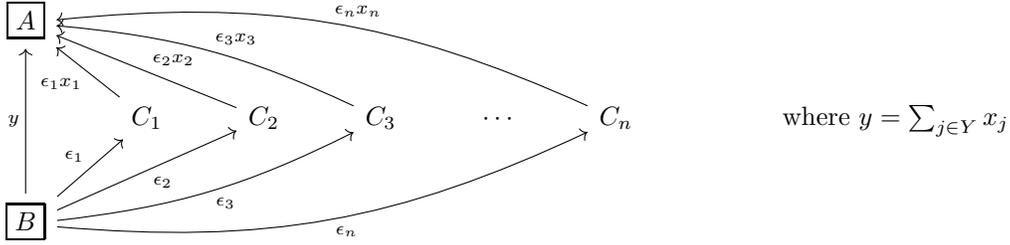
\vspace{5mm}

\begin{center}

\begin{tikzcd}
\boxed{A} \\
& C_1 \ar[ul, "\epsilon_1 x_1"] & C_2 \ar[ull, "\epsilon_2 x_2"'] & C_3 \ar[ulll, "\epsilon_3 x_3"', bend right = 10]& \cdots &  C_n \ar[ulllll,  "\epsilon_n x_n"', bend right = 15] && \text{where } y = \sum_{j \in Y} x_j \\
\boxed{B} \ar[ur, "\epsilon_1"]  \ar[urr, "\epsilon_2"'] \ar[urrr, "\epsilon_3"', bend right = 10] \ar[urrrrr,  "\epsilon_n"', bend right = 15] \ar[uu, "y"]
\end{tikzcd}

\end{center}

\vspace{5mm}

That means that if $k \notin \{0, 1\} \cup X$, the only way for $\mu(Q)$ to contain an arrow with weight $k$ is for $k$ to be the weight of the arrow between $B$ and $A$. That means that $k$ is present in some quiver equivalent to $Q$ if and only if $k$ is a subset-sum of $X$. The result follows from \textsf{NP-hardness} of \textsc{Subset Sum}. \qed

\subsection*{Proof of Theorem \ref{t:strongnphard}}

We use the following formulation of the \textsc{3-partition} problem:

\smallskip

\begin{tabular}{l l}
\multicolumn{2}{l}{\textsc{\large{3-Partition } } } \\
\textbf{Input:}  & $n \geq 3$, and $\mathcal{X} \subseteq \binom{[n]}{3}$ .							\\ 
\textbf{Decide:} \quad & $\exists \mathcal{A} \subseteq \mathcal{X}$ ~ such that ~ every $i \in [n]$ is contained in exactly one $A \in \mathcal{A}$?

\end{tabular}

\smallskip

 given a positive integer $n$ and a subset $\mathcal{X} \subseteq \binom{[n]}{3}$, does there exist a partition of $[n]$ into elements of $\mathcal{X}$? This is \textsf{strongly NP-hard} (see e.g. [GJ], \S A3.1).  Without loss of generality, we may assume that each element of $[n]$ is in at least one of the elements of $\mathcal{X}$.

Given $n$ and $X$, we construct a quiver with vertices:
\[
\underbrace{A_1, \dots, A_n}_{\text{frozen}}, \{B_X\}_{X \in \mathcal{X}}, \underbrace{C}_{\text{frozen}}
\]
Take the following edges:
\begin{itemize}
\item[] One edge from $A_i$ to $B_X$ if $i \in X$.
\item[] One edge from $B_X$ to $C$ for each $X$.
\item[] One edge from $C$ to $A_i$ for each $i$.
\end{itemize}
The resulting quiver has this shape:

\smallskip
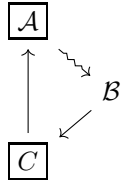
\begin{center}

\begin{tikzcd}[sep=small]
\boxed{\mathcal{A}} \ar[dr, rightsquigarrow]\\
& \mathcal{B}  \ar[dl]\\
\boxed{C} \ar[uu]
\end{tikzcd}

\end{center}

\smallskip
Here $\mathcal{A}$ and $\mathcal{B}$ represent the sets of vertices of the form $A_i$ and $B_X$ respectively. There is only one vertex labeled $C$. The solid arrows represent one arrow between every pair of vertices from the respective sets. And the squiggly arrow between  $\mathcal{A}$ and $\mathcal{B}$ represents arrows between an $A_i$ and a $B_X$ if and only if $i \in X$.

If a partition $\mathcal{P} \subset \mathcal{X}$ exists, we can apply the mutations $\{\mu_{B_P} : P\in \mathcal{P} \}$ followed  which will eliminate all icebound edges.

More generally, note that all mutations commute. Moreover, they are all involutions. So we need only consider the effect of using mutations at most once. In that case, we eliminate the icebound edges if and only if the mutations we use correspond to a partition of $[n]$. We have therefore reduced the problem to \textsc{3-partition}. Because \textsc{3-partition} is \textsf{strong NP-hard}, the result follows. \qed

\smallskip

\subsection*{Proof of Theorem \ref{t:alternate}}

First, we note that it suffices to prove the case where $Q$ has exactly four vertices. This is because, for any subset $Q' \subset V(Q)$ of size 4 containing both $C$ and $D$, the action of $\mu_C$ and $\mu_D$ commutes with restriction to $Q'$. 

Let the other two vertices in $Q'$ be $A$ and $B$. It also suffices to consider the case where $A$ and $B$ start with no arrows between them. If $C$ and $D$ have no arrows between them to start, then the statement is trivial. If $C$ and $D$ have one arrow between them, then it is an easy computation to check that $(\mu_D\mu_C)^{10}Q = Q$.

So assume there are $\alpha \geq 2$ arrows from $C$ to $D$. One possible case consists of arrows from $A$ to $C$ and from $D$ to $B$
 Then we can write down the first few quivers that we get:

\smallskip

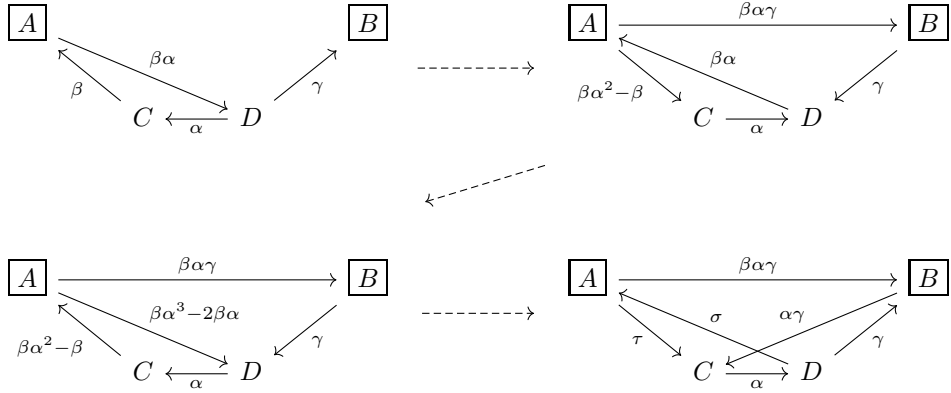
\begin{center}
\begin{tikzcd}
\boxed{A} \ar[drr, "\beta\alpha"]& & & \boxed{B} 
&& \boxed{A} \ar[rrr, "\beta\alpha \gamma"] \ar[dr, "\beta\alpha^2 - \beta"'] & & & \boxed{B}  \ar[dl, "\gamma"] \\
& C \ar[ul, "\beta"] & \ar[l, "\alpha"] D  \ar[ur, "\gamma"']& \ar[rr, dashed, shift left = 8, shorten=5mm]
&\ &\ \ar[dll, dashed, shorten = 5mm, shift left=5]& C \ar[r, "\alpha"']& D \ar[ull, "\beta\alpha"']& \\ \ &\ &\ &\ \ar[rr, dashed, shift right = 20, shorten=5mm] &\ &\ &\ &\ \\
\boxed{A} \ar[rrr, "\beta\alpha \gamma"]  \ar[drr, "\beta\alpha^3 - 2\beta\alpha"]& & & \boxed{B}  \ar[dl, "\gamma"]
&&\boxed{A} \ar[rrr, "\beta\alpha \gamma"]  \ar[dr, "\tau"']& & & \boxed{B}   \ar[dll, "\alpha \gamma"']\\
& C \ar[ul, "\beta\alpha^2 - \beta"] & D  \ar[l, "\alpha"]&
&&& C  \ar[r, "\alpha"'] & D \ar[ur, "\gamma"']\ar[ull, "\sigma"'] &
\end{tikzcd}

\end{center}

where $\sigma = \beta\alpha^3 - 2\beta\alpha$ and $\tau = \beta\alpha^4 - 2\beta\alpha^2 - \beta\alpha^3 + \beta$. Note that these are both positive since $\alpha \geq 2$.

Consider the quivers $Q_1(x, y, z, w)$ and $Q_2(p, q, r, s)$:

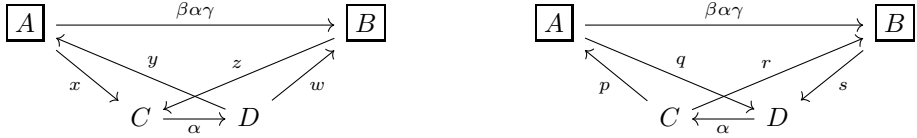
\smallskip

\begin{center}

\begin{tikzcd}
\boxed{A} \ar[rrr, "\beta\alpha \gamma"] \ar[dr, "x"'] & & & \boxed{B} \ar[dll, "z"'] 
&&\boxed{A} \ar[rrr, "\beta\alpha \gamma"] \ar[drr, "q"]& & & \boxed{B} \ar[dl, "s"]\\
& C \ar[r, "\alpha"']& D \ar[ur, "w"'] \ar[ull, "y"']&
&&& C \ar[ul, "p"] \ar[urr, "r"]& \ar[l, "\alpha"] D &
\end{tikzcd}

\end{center}

\smallskip

We claim that all future quivers will be of one of these two forms and thus $\beta\alpha \gamma$ is the only thing that appears on top. We can compute that 
\[
\mu_C \big(Q_1(x, y, z, w)\big) = Q_2(x, \alpha x - y, z, \alpha z - w)
\]
so long as $\alpha x > y$ and $\alpha z > w$. Next we apply $\mu_D$ to find

\begin{align*}
\mu_D\Big(\mu_C\big(Q_1(x, y, z, w)\big)\Big) &= Q_1\big(\alpha(\alpha x - y) - x ,\alpha x - y, \alpha(\alpha z - w) - z, \alpha z - w  \big) \\
& : = Q_1(x', y', z', w')
\end{align*}

this time assuming $\alpha(\alpha x - y) > x$ and $\alpha(\alpha z - w) > w$. This is a stronger condition than the previous. Note that our conditions are equivalent to
\[
\frac{\alpha}{\alpha^2 - 1} < \min \left( \frac{x}{y}, \frac{z}{w} \right)
\]
which is satisfied by our original picture.

However, we have computed

\[
\frac{x'}{y'} = \frac{\alpha(\alpha x - y) - x}{\alpha x - y} 
= \alpha - \frac{x}{\alpha x - y} = \alpha - \frac{x/y}{\alpha(x/y) - 1}.
\]
So the problem reduces to iteratively applying the function
\[
f(t) = \alpha - \frac{t}{\alpha t - 1}
\]
and it is easy to see that this converges to a limit of 
\[
\frac{x}{y} = t = \frac{1}{2}\left(\alpha - \sqrt{\alpha^2 - 4}\right)
\]
A similar picture holds for the other three starting positions. \qed
\smallskip
\section{Final Remarks and Open Problems} \label{s:finalremarks}

\subsection{Undecidability.}  The paper \cite{FN23+} does show the existence of small quivers which are nonetheless polynomially far apart with respect to mutation. Of course, undecidability is far stronger. Suppose, for example, that it is undecidable whether or not two quivers are mutation equivalent. Then, there would exist quivers $Q_1$ and $Q_2$ with $a_1$ and $a_2$ arrows respectively such that the shortest sequence of mutations taking one to the other has length
\[
\ell \geq  (\Tow (\Tow(a_1 + a_2 ))
\]
where $\Tow(k)$ is a tower of 2's of length $k$. Put fancifully, this means there is no limit as to how far into the sky one has to go in order to show that two quivers are mutation equivalent. Note that Theorem~\ref{t:alternate} shows that more than two mutable vertices are needed for this to happen.

\subsection{Knots and Plabic Graphs.} Deep connections exist between quiver mutation equivalence and knot theory including via \emph{plabic graphs} (see e.g. \cite{GL}, \cite{BS}, \cite{STWZ} or \cite{FPST22}). For knots and links, upper bounds exist for the number of Reidemeister moves needed to show equivalence.

Here is the best known bound due to \cite{CL14}. Suppose $D_1$ and $D_2$ are diagrams of the same link or knot. Let their crossing numbers be $c_1$ and $c_2$ respectively. Then there exists a sequence of Reidemeister moves taking $D_1$ to $D_2$ of length at most
\[
\Tow(C)^{(c_1 + c_2)} \qquad \text{where } C = \left( 10^{10^6}\right)^{(c_1 + c_2)}
\]
Again, $\Tow(k)$ is a tower of 2s of length $k$. It may well be the case that such a bound exists for mutation equivalence of quivers as well.

\subsection{Quiver invariants} \label{ss:pnp} Fix a quiver $Q$ and some $k  \in \mathbb{N}$. Fomin's question in the introduction asks for an algorithm to determine whether there exists a quiver $Q' \in [Q]$ such that two vertices in $Q'$ have exactly $k$ arrows between them.

One hope is that determinantal invariants would be able to answer these questions. Our results suggest that one should investigate the quivers used in the construction of Theorems \ref{t:nphard} and \ref{t:strongnphard}.

\subsection{Mutable and immutable vertices.} Frozen vertices and arrows between them are essential to the proofs of Theorems \ref{t:nphard} and \ref{t:strongnphard}. It would be interesting to see whether the number of frozen vertices can be reduced.

\subsection{Other properties of quivers}. There are many other questions about quivers for which an algorithmic test would be of interest. For instance, one could ask whether a quiver is \emph{mutation-acyclic}, that is, mutation equivalent to an acylic quiver. Much work remains to be done in this area.

\subsection{Quiver gadgets} \label{ss:gadgets} Embedding difficult problems into quiver mutation equivalence requires the construction of quivers whose mutations can be controlled. Many questions even about simple quivers remain unanswered. In particular, one method would be to embed \emph{Hilbert's 10th problem} or the \emph{Post correspondence problem} into quivers (see e.g. \cite{PS}). We are still far away from this.

To illustrate, we give a natural possible generalization of Theorem  \ref{t:alternate}. Let $Q$ be a quiver of the following form:

\begin{center}
\smallskip

\begin{tikzcd}
\boxed{A} \ar[r, "x_0"] & C_1 \ar[r, "x_1" ] & \cdots \ar[r, "x_{k - 1}"] & C_k \ar[r, "x_k"]& \boxed{B}
\end{tikzcd}

\smallskip
\end{center}

Then we conjecture that for all $Q'$ which is mutation equivalent to $Q$, the number of arrows between $A$ and $B$ is always $0$ or $x_0x_1\cdots x_k$. The cases $k = 0$ and $k = 1$ are trivial, and the case $k = 2$ is proven in Theorem \ref{t:alternate}. However, the general case is open.

\smallskip

\section*{Acknowledgments}
 I would like to thank Sergey Fomin and Bernhard Keller for helpful suggestions. .I am grateful to Pavel Galashin and Nikita Gladkov for fruitful conversations. I would especially like to thank my advisor Igor Pak for helpful comments and guidance.

\end{document}